\documentstyle[12pt]{article}
  \newcommand{\const}{\rm const}
  
  \newcommand{\supp}{\rm supp}

\begin{document}

   \begin{center}

  \   {\bf   Fourier transform acting on the functions }\\

\vspace{4mm}

   \  {\bf defined in the infinite LCA groups.  }\\

\vspace{4mm}

 \ {\bf   A   Grand Lebesgue Spaces approach. }\\

\vspace{4mm}

  \   {\bf Ostrovsky E., Sirota L.}\\

\vspace{4mm}

 Israel,  Bar - Ilan University, department  of Mathematic and Statistics, 59200, \\

\vspace{4mm}

E - mails: eugostrovsky@list.ru, \hspace{5mm} sirota3@bezeqint.net \\

\vspace{5mm}

  {\bf Abstract} \\

\vspace{4mm}

 \end{center}

 \  We derive  in this article the  exact norm in the Grand Lebesgue Spaces (GLS) estimates for Fourier transform  acting on the
functions defined in the infinite local compact Abelian (LCA) group, compact  or discrete.  \\

\vspace{4mm}


\vspace{5mm}

\section{1. \ Definitions.  Notations. Previous results.  Statement of problem.}
 \vspace{3mm}

 \ Let  $  \ X \ $ be infinite  local compact Abelian (LCA) group,  compact or discrete; $  \ T \ $ be unit circle on the plane,

$$
 Y   = \{ \gamma: X \to T  \} \eqno(1.0)
$$
be its dual group,  $  \ \alpha(\cdot) \ $ be the Haar's measure on the group $  \ X; \  \ \beta(\cdot) \ $  be the Haar's measure
on the dual group $ \ Y; \  F[f] (\cdot) = \hat{f}(\cdot)  \ $  be the Fourier transform of the measurable function $ \ f: X \to C: $

$$
F[f] (\gamma) = \hat{f}(\gamma)  \stackrel{def}{=} \int_X f(x) \ \gamma(-x) \ \alpha(dx).  \eqno(1.1)
$$
 \ The reciprocal norming of the Haar's measures  $ \ \alpha, \ \beta  $   will be presumed such that

$$
f(x) = \int_Y \hat{f}(\gamma) \ \gamma(x) \ \beta(d \gamma).   \eqno(1.2)
$$

 \ Put  also

$$
A := \alpha(X), \hspace{5mm} B := \beta(Y).
$$

 \ Recall that the classical Lebesgue - Riesz $ \ L_p = L_p(X), \ p \ge 1 \ $ norm for the (measurable) function $ f: X \to C $  is defined as ordinary
by the formula

$$
|f|L_p(X) = |f|_p :=  \left[ \int_X  |f(x)|^p \ \alpha(dx)  \right]^{1/p}.
$$
 \ The correspondent $ L_q = L_q(Y), \ q \ge 1 $ norm  for the Fourier transform of the function $  f(\cdot) $ has  by definition a form

$$
 |\hat{f}|_q(Y) =   |\hat{f}|_q :=  \left[ \int_Y  | \hat{f}(\gamma)|^q \ \beta(d  \gamma)   \right]^{1/q}.
$$

  \ Denote for each number $  p, \ p > 1 $ by $  p' $ its conjugate number:

$$
p' = \frac{p}{p-1}, \hspace{4mm} \Leftrightarrow \frac{1}{p}  + \frac{1}{p'} = 1.
$$
 \ Analogously $  \ q' = q/(q-1), \ q > 1.  $\par

 \ Denote  following the authors of the very interest article  [7]  by $ \ K(p,q) \ $ the norm of the Fourier transform in the
$  \ L_q(Y) \to L_p(X) $ sense:

$$
K(p,q) := \sup_{ |f|_p = 1  } |\hat{f}|_q  = \sup_{  0 \ne f \in L_p(X) }  \frac{|\hat{f}|_q(Y)}{|f|_p(X)}=
$$

$$
||F[\cdot]|| \left[L_p(X) \to L_q(Y) \right],  \eqno(1.3)
$$
so that

$$
 |\hat{f}|_q(Y)  \le K(p,q) \ |f|_p(X).
$$

 \ Let us introduce the following domain

$$
Q := \left\{ (p,q) \in (0,\infty)^2: \ 1/p + 1/q \le 1, \ q \ge 2 \right\}.  \eqno(1.4)
$$

 \ Mokshay Madiman and Peng Xu  in the article [7]  calculated the exact value of the variable  $ \ K(p,q). \ $
Namely, they proved that  if $  (p,q) \in Q, $ then

$$
K(p,q) = A^{1 - 1/p - 1/q }   \eqno(1.5)
$$
and $ \ K(p,q) = \infty \ $ otherwise.  \par

 \ The  ''conjugate''   proposition has the following form. Define the following domain and variable:

$$
\hat{Q} := \{ (p,q): p \in (0,2), \ 1/p + 1/q \ge 1, \ p \le 2 \};
$$

$$
\hat{K}(q,p) := B^{ 1/p + 1/q -1 },  \ (q,p) \in \hat{Q}; \  \hat{K}(q,p)  = \infty,  \ (q,p) \notin \hat{Q};  \eqno(1.6)
$$
then

$$
\hat{K}(q,p) = \sup_{ 0 \ne f \in L_p(X)  } \frac{| f|_p(X)}{|\hat{f}|_q(Y)},  \eqno(1.7)
$$
so that

$$
|f|_p(X) \le B^{1/p + 1/q - 1} \ |\hat{f}|_q(Y), \  (p,q) \in \hat{Q},  \eqno(1.8)
$$
and the ``constant''  $  B^{1/p + 1/q - 1} $ is the best possible. \par

\vspace{4mm}

 \ {\bf   Our target in this article  is to extend the estimates (1.5)  - (1.8)  into the case when instead the classical Lebesgue - Riesz
spaces $  \ L_p, \ L_q  $ stands the so - called Grand Lebesgue Spaces (GLS). }\par

\vspace{3mm}

\section{2. \ Grand Lebesgue Spaces (GLS). Fundamental functions. }

 \vspace{3mm}

 \ Let $  Z = (Z, M, \mu) $ be measurable space with non - trivial sigma - finite measure $ \ \mu. \  $
 \ Let also  $  \psi = \psi(p), \ p \in [1, b), \ b = \const \in (1,\infty]  $ (or   $ p \in [1,b] $ ) be certain bounded
from below:  $  \ \inf \psi(p)  > 0 $ continuous inside the  {\it  semi - open} interval $   \ p \in [1, b) $ numerical function. We can and will suppose
$   \ b = \sup \{p, \psi(p) < \infty\},  \ $ so that  $ \  \supp \ \psi = [1, b) \  $  or $ \ \supp \ \psi = [1, b]. \ $ The set of all such a functions will be denoted by
$ \ \Psi(b) = \{  \psi(\cdot)  \}; \ \Psi := \Psi(\infty).  $\par

 \ By definition, the (Banach) Grand Lebesgue Space  \ (GLS)  \ space   $  \ G\psi = G\psi(b)  $ consists on all the   complex numerical  valued
measurable functions $ \zeta  $ defined on our measurable space  and having a finite norm

$$
||\zeta|| = ||\zeta||G\psi \stackrel{def}{=} \sup_{p \in [1,b)} \left\{ \frac{|\zeta|_p}{\psi(p)} \right\}, \eqno(2.0)
$$
here

$$
|\zeta|_p = |\zeta|L_p(Z).
$$

 \ These spaces  are Banach functional space, are complete, and rearrangement invariant in the classical sense, see [1],  chapters 1, 2;
and were  investigated in particular in  many  works, see e.g.  [3],  [4], [5], [6],  [8],  [15], [16], [17].
 We refer here  some  used in the sequel facts about these spaces  and supplement more. \par

\vspace{4mm}

 \  Suppose temporarily that the measure $  \ \mu \ $ is probabilistic:  $ \  \mu(A) = {\bf P}(A); \  \mu(Z) = {\bf P}(Z) =1. \ $

\vspace{3mm}

\  It is known that if  $  \ \zeta \ne 0,  $

  $$
  {\bf P}(|\zeta| > y) \le \exp \left( \ - v_{\psi}^* (\ln (y/||\zeta||) \ \right), \ y \ge e \cdot ||\zeta||. \eqno(2.2)
  $$

  \ Conversely,  the last inequality may be reversed  in the following version: if

$$
{\bf P}(|\zeta| > y) \le \exp \left(-v_{\psi}^* (\ln (y/K) \right). \ y \ge e \cdot K, \ K = \const \in (0,\infty),
$$
and if the function $ v_{\psi}(p), \ 1 \le p < \infty  \ $  is positive, continuous, convex and such that

$$
\lim_{p \to \infty} \ln \psi(p) = \infty,
$$
then  $ \zeta \in G\psi $ and besides $  ||\zeta||  \le C(\psi) \cdot K:  $

$$
||\zeta||G\psi \le C_1 ||\zeta||L(M)  \le C_2 ||\zeta||G\psi, \ 0 < C_1 < C_2 < \infty. \eqno(2.3)
$$

 \  Furthermore,  let now $  \eta = \eta(z), \ z \in W $ be arbitrary family  of random variables  defined on any set $ W  $ such that

$$
\exists b \in (1,\infty] \ \forall p \in [1,b)  \ \Rightarrow  \psi_W(p) := \sup_{z \in W} |\eta(z)|_p  < \infty. \eqno(2.4)
$$
 \ The function $  p \to \psi_W(p)  $ is named as a {\it  natural} function for the  family  of random variables $  W.  $  Obviously,

$$
\sup_{z \in W} ||\eta(z)||G\Psi_W = 1.
$$

\vspace{4mm}

 \ {\bf   Definition 2.1. }  The {\it fundamental function} for GLS $ \ G\psi_b \hspace{4mm} \phi[G\psi](\delta), \ \delta \in (0,\infty) $ is defined
by a formula

$$
\phi[G\psi](\delta)  := \sup_{p \in [1,b)} \left\{  \frac{\delta^{1/p}}{ \psi(p) }  \right\}. \eqno(2.5)
$$
 \ This notion play a very important role in the Functional Analysis, theory of Fourier series, Operator Theory, Theory of Random Processes
 etc., see  the classical monograph [1]. For the GLS  this function was investigated in the preprint  [17].\par

 \ {\bf   Definition 2.2. }  The {\it low truncated fundamental function} for the GLS $ \ G\psi_b, $ \hspace{3mm} namely,
$ \ \phi_s[G\psi](\delta), \ \delta \in (0,\infty), 0 < s < b \ $  is defined  by a formula

$$
\phi_s[G\psi](\delta)  := \sup_{p \in [s,b)} \left\{  \frac{\delta^{1/p}}{ \psi(p) }  \right\}, \ 1 \le s < b. \eqno(2.5a)
$$

 \ {\bf   Definition 2.3. }  The {\it upper truncated fundamental function} for the GLS
$ \ G\psi_b, $ \hspace{4mm}  indeed: $ \ \phi^s[G\psi](\delta), \ \delta \in (0,\infty), 0 < s < b \ $ is defined  by a formula

$$
\phi^s[G\psi](\delta)  := \sup_{p \in [s,b)} \left\{  \frac{\delta^{1/p}}{ \psi(p) }  \right\}, 1 \le s < b. \eqno(2.5b)
$$

\vspace{4mm}

\section{ 3. Main results.}

 \vspace{3mm}

 \ {\bf A. Compact case.}  \par

 \ Suppose for beginning  that the LCA group $  \  X \ $ is compact. Then $  A = \alpha(X) \  < \infty $ and the dual group is discrete with correspondent
 counting Haar's measure $  \ \beta. \ $\par

 \  Assume in this sub - section that the measurable  function $  f: X \to C $ belongs to some $ \ G\psi_b \ $ space:

$$
|f|_p \le \psi(p) \cdot ||f||G\psi_b, \ b = \const \in (1, \infty]; \ p \in [1, b). \eqno(3.0)
$$
 \ We can and will suppose without loss of generality $   |f||G\psi_b = 1. $ Of course, the function $   \psi = \psi(p) $ may be picked as a
natural function fort the function $ f: \ \psi(p) = |f|_p(X), \ $ if it is non - trivial.\par

 \ Let $ \ (p,q) \ $ be  the numbers  from the set  $  Q; $ this implies

$$
q \ge 2, \ q \ge p'; \hspace{5mm} \Leftrightarrow   q \ge \max(2,p'). \eqno(3.1)
$$

 \ It follows  from the result of Mokshay Madiman and Peng Xu  [7],   \hspace{4mm}  (1.5)

$$
|\hat{f}|_q \le A^{1 - 1/p - 1/q }  \cdot |f|_p \le A^{1 - 1/p - 1/q } \ \psi(p), \eqno(3.2a)
$$
or equally

$$
A^{1/q - 1} |\hat{f}|_q \le A^{-1/p} \ \psi(p). \eqno(3.2b)
$$
 \ We deduce from (3.2a)  taking maximum over  $  p  $ and denoting $ \ t(q)  = \min(q',2): $

$$
A^{1/q - 1} |\hat{f}|_q \le  \inf_{p \le t(q)}  \left[ A^{-1/p} \ \psi(p) \right]. \eqno(3.3)
$$

 \ The right - hand side of (3.3) may be rewritten as follows:

$$
  \inf_{p \le t(q)} \left[ A^{-1/p} \ \psi(p) \right] =  \sup_{p \le t(q)} \left\{ \frac{A^{1/p}}{\psi(p)}   \right\}^{-1} =
$$

$$
\left\{ \phi[G\psi]_{t(q)} (A) \right\}^{-1},
$$
therefore

$$
A^{-1}  |\hat{f}|_q  \le \frac{A^{1/q} }{\phi[G\psi]_{t(q)} (A) }. \eqno(3.4)
$$

 \ Suppose in addition that the function $ \phi[G\psi]_{t(q)} (A)  $  allows a following representation (factorization)

$$
\phi[G\psi]_{t(q)} (A) = \frac{\theta_A(q)}{\nu_A(q)}, \ q \ge 2, \eqno(3.5)
$$
where both the functions  $ \theta_A(\cdot), \  \nu_A(\cdot) \ $ are from the set $ \ G[2, d); \ d = \const \in (2, b). $ \par

 \ Then the inequality (3.4) may be rewritten as follows.

$$
\frac{ |\hat{f}|_q}{\nu_A(q)}  \le  A \cdot \frac{A^{1/q}}{\theta_A(q)}, \eqno(3.6)
$$
and we obtain  after taking the maximum over both the sides  of the inequality (3.6)

$$
||\hat{f}||G\nu_A \le A \cdot \phi[G\theta_A](A). \eqno(3.7)
$$

 \ Thus, we obtained really the following result. \par

\vspace{4mm}

{\bf Theorem 2.1.} We conclude under formulated in this section conditions, restrictions and notations, in particular,
under the condition (3.5),

 $$
||\hat{f}||G\nu_A \le A \cdot \phi[G\theta_A](A) \cdot ||f||G\psi. \eqno(3.8)
$$

\vspace{4mm}

 \ {\bf B. Discrete case.}  \par

 \ Suppose  now that the LCA group $  \  X \ $ is discrete. Then $  B = \beta(Y) \  < \infty $ and the dual group $ \ Y \  $ is compact with
correspondent  Haar's measure $  \ \alpha. \ $\par

  \ Our considerations in this subsection are alike the foregoing ones  in the investigation of theorem 2.1. \par

 \  Assume  here that the measurable  function $ \hat{f}: \ Y \to C $ belongs to some $ \ G\psi_b \ $ space:

$$
|\hat{f}|_q \le \psi(q) \cdot ||f||G\psi_b, \ b = \const \in (2, \infty]; \ q \in [2, b), . \eqno(3.9)
$$
 \ We can and will suppose without loss of generality $   ||\hat{f}||G\psi_b = 1. $ Of course, the function $   \psi = \psi(q) $ may be choosed as a
natural function fort the function $ f: \ \psi(q) = |\hat{f}|_q(Y), \ $ if  the last function is non - trivial.\par

 \ Let $ \ (p,q) \ $ be  the numbers  from the set  $  \hat{Q}; $ this implies

$$
1 \le p \le 2, \ p \le q'; \hspace{5mm} \Leftrightarrow   1 \le p \le \min(2,q'). \eqno(3.10)
$$

 \ It follows again  from the second result of Mokshay Madiman and Peng Xu [7], theorem 4.1;   \hspace{4mm}  (1.7)

$$
|f|_p \le B^{-1 + 1/p + 1/q }  \cdot | \hat{f}|_q \le B^{-1 + 1/p + 1/q } \ \psi(q), \eqno(3.11)
$$
or equally

$$
B^{ 1 - 1/p} |f|_p \le B^{1/q} \ \psi(q). \eqno(3.11a)
$$
 \ We deduce from (3.11a)  taking maximum over  $  p  $ and denoting $ \ s(p)  = \max(p',2): $

$$
B^{1 - 1/p } |f|_p \le  \inf_{q \ge s(p)}  \left[ B^{1/q} \ \psi(q) \right]. \eqno(3.12)
$$

 \ The right - hand side of (3.12) may be rewritten as follows:

$$
  \inf_{q \ge s(p)} \left[ B^{1/q} \ \psi(q) \right] =  \sup_{q \ge s(p)} \left\{ \frac{B^{-1/q}}{\psi(q)}   \right\}^{-1} =
$$

$$
\left\{ \phi[G\psi]^{s(p)} \left( B^{-1} \right) \right\}^{-1},
$$
therefore

$$
B^{-1}  |f|_p  \le \frac{B^{1/p} }{\phi[G\psi]^{s(p)} (B^{-1}) }. \eqno(3.13)
$$

 \ Suppose in addition that the function $ \phi[G\psi]^{s(p)} (B)  $  allows a following representation (factorization)

$$
\phi[G\psi]^{s(p)} (B^{-1}) = \frac{\tau_B(q)}{\kappa_B(p)}, \ p \ge q', \eqno(3.14)
$$
where both the functions  $ \tau_B(\cdot), \  \kappa_B(\cdot) \ $ are from the set $ \ G[1, b); \ b = \const \in (1,2]. $ \par

 \ Then the inequality (3.13) may be rewritten as follows.

$$
\frac{ |f|_p}{\kappa_B(p)}  \le  B^{-1} \cdot \frac{B^{-1/p}}{\tau_B(p)}, \eqno(3.15)
$$
and we obtain  after taking the maximum over both the sides  of the inequality (3.15)

$$
||f||G\kappa_B \le B^{-1} \cdot \phi[G\tau_B] \left(B^{-1} \right). \eqno(3.16)
$$

 \ Thus, we obtained really the following result. \par

\vspace{4mm}

{\bf Theorem 2.2.} We conclude under formulated in this section conditions, restrictions and notations, in particular,
under the condition (3.14),

 $$
||f||G\kappa_B \le B^{-1} \cdot \phi[G\tau_B] \left(B^{-1} \right) \cdot ||f||G\psi. \eqno(3.17)
$$

\vspace{4mm}

 \ {\bf Remark 3.1.} It is easily to verify that  the source estimated of Mokshay Madiman and Peng Xu [7] are the
particular cases of obtained here (3.8) and  (3.17), as long as the classical Lebesgue - Riesz spaces $  \ L_p \ $ are the particular
cases of Grand Lebesgue Spaces.  \par
 \ This implies as a consequence that both the constants in the right - hands sides of inequalities
(3.8) and  (3.17) are in general case non - improvable. \par

\vspace{4mm}

\section{ 4. \ Concluding remarks.}

 \vspace{3mm}

 \ It is interest by our opinion to generalize obtained in \cite{Mokshay1} and in
this report estimates to the non - commutative infinite local compact groups,
as well as in the other rearrangement invariant  spaces builded over $  X $ or $ Y. $ \par

 \vspace{6mm}

 {\bf References.}

 \vspace{4mm}

{\bf 1. Bennet C., Sharpley R.}  {\it  Interpolation of operators.} Orlando, Academic
Press Inc., (1988). \\

 \vspace{3mm}

{\bf 2.  Buldygin V.V., Kozachenko Yu.V. }  {\it Metric Characterization of Random
Variables and Random Processes.} 1998, Translations of Mathematics Monograph, AMS, v.188. \\

 \vspace{3mm}

 {\bf 3. A. Fiorenza.}   {\it Duality and reflexivity in grand Lebesgue spaces. } Collect. Math.
{\bf 51,}  (2000), 131  - 148. \\

 \vspace{3mm}

{\bf  4. A. Fiorenza and G.E. Karadzhov.} {\it Grand and small Lebesgue spaces and
their analogs.} Consiglio Nationale Delle Ricerche, Instituto per le Applicazioni
del Calcoto Mauro Picone”, Sezione di Napoli, Rapporto tecnico 272/03, (2005).\\

 \vspace{3mm}

{\bf 5.  T. Iwaniec and C. Sbordone.} {\it On the integrability of the Jacobian under minimal
hypotheses. } Arch. Rat.Mech. Anal., 119, (1992), 129-143. \\

 \vspace{3mm}

{\bf 6. Kozachenko Yu. V., Ostrovsky E.I. }  (1985). {\it The Banach Spaces of random Variables of sub-gaussian Type. } Theory of Probab.
and Math. Stat. (in Russian). Kiev, KSU, 32, 43-57. \\

 \vspace{3mm}

{\bf 7.  Mokshay Madiman and Peng Xu.} {\it  The norm of the Fourier transform on compact or discrete abelian groups. }
arXiv:1611.04692 v1 [math.CA] 15 Nov 2016. \\

 \vspace{3mm}

{\bf 8. Ostrovsky E.I. } (1999). {\it Exponential estimations for Random Fields and its
applications,} (in Russian). Moscow-Obninsk, OINPE. \\

 \vspace{3mm}

{\bf 9. Ostrovsky E. and Sirota L.} {\it Vector rearrangement invariant Banach spaces
of random variables with exponential decreasing tails of distributions.} \\
 arXiv:1510.04182v1 [math.PR] 14 Oct 2015 \\

 \vspace{3mm}

{\bf 10. Ostrovsky E. and Sirota L.}  {\it Non-asymptotical sharp exponential estimates
for maximum distribution of discontinuous random fields. } \\
 arXiv:1510.08945v1 [math.PR] 30 Oct 2015 \\

 \vspace{3mm}

{\bf 11. Ostrovsky E.I.}  {\it About supports of probability measures in separable Banach
spaces.} Soviet Math., Doklady, (1980), V. 255, N 6, p. 836-838, (in Russian).\\

 \vspace{3mm}

 {\bf 12. Ostrovsky E. and Sirota L.} {\it Criterion for convergence almost everywhere, with applications.} \\
arXiv:1507.04020v1 [math.FA] 14 Jul 2015. \\

 \vspace{3mm}

{\bf 13. Ostrovsky E. and Sirota L.}  {\it Schl\"omilch and Bell series for Bessels functions, with probabilistic applications.} \\
 arXiv:0804.0089v1 [math.CV] 1 Apr 2008 \\

 \vspace{3mm}

{\bf 14. Ostrovsky E. and Sirota L. } {\it Sharp moment estimates for polynomial martingales. } \\
arXiv:1410.0739v1 [math.PR] 3 Oct 2014 \\

 \vspace{3mm}

{\bf 15. Ostrovsky E., Rogover E. } {\it Exact exponential bounds for the random field
maximum distribution via the majorizing measures (generic chaining).} \\
 arXiv:0802.0349v1 [math.PR] 4 Feb 2008 \\

 \vspace{3mm}

{\bf 16. Ostrovsky E. and Sirota L. } {\it   Entropy and Grand Lebesgue Spaces approach for the problem  of Prokhorov - Skorokhod
continuity of discontinuous random fields. }\\
arXiv:1512.01909v1 [math.Pr] 7 Dec 2015 \\

 \vspace{3mm}

{\bf 17. Ostrovsky E. and Sirota L. } {\it  Fundamental function for Grand Lebesgue Spaces.  }
arXiv:1509.03644v1  [math.FA]  11 Sep 2015 \\

\end{document}